\magnification=\magstep1

\font\bfw=cmbx12
\hoffset=.4truecm
\voffset= 1truecm
\hsize=16.5truecm
\vsize=24truecm
\leftskip=0.5cm
\baselineskip=12pt
\hfuzz=18pt
\parindent=1truecm
\parskip=0.2truecm
\parskip=0.2truecm
\input psfig


\def\be{\begin{equation}}
\def\ee{\end{equation}}
\def\bea{\begin{eqnarray}}
\def\eea{\end{eqnarray}}

\def\hfl#1#2{\smash{\mathop{\hbox to 12mm{\rightarrowfill}}
\limits^{\scriptstyle#1}_{\scriptstyle#2}}}
\def\vfl#1#2{\llap{$\scriptstyle #1$}\left\downarrow
\vbox to 6mm{}\right.\rlap{$\scriptstyle#2$}}
\def\ihfl#1#2{\smash{\mathop{\hbox to 12mm{\nearrowfill}}
\limits^{\scriptstyle#1}_{\scriptstyle#2}}}
\def\ibfl#1#2{\smash{\mathop{\hbox to 12mm{\searrowfill}}
\limits^{\scriptstyle#1}_{\scriptstyle#2}}}
\def\diagram#1{\def\normalbaselines{\baselineskip=0pt
\lineskip=10pt\lineskiplimit=1pt}  \matrix{#1}}
\def\adots{\mathinner{\mkern2mu\raise1pt\hbox{.}
\mkern3mu\raise4pt\hbox{.}\mkern1mu\raise7pt\hbox{.}}}
\def\sqr#1#2{{\vcenter{\vbox{\hrule height.#2pt
        \hbox{\vrule width.#2pt height#1pt \kern#1pt
           \vrule width.#2pt}
        \hrule height.#2pt}}}}

\baselineskip=14pt
\line{\hfill                                                     q-alg/9509016}
\line{\hfill                                                   CPTH-RR361.0695}
\vskip 4truecm
\centerline{\bfw ATYPICAL REPRESENTATIONS OF ${\cal U}_{\displaystyle{q}}
(sl(N))$ AT}
\centerline{\bfw ROOTS Of UNITY}
\smallskip
\centerline{Boucif Abdesselam $^{}$\footnote{$^{1}$}
{\it E-mail adress: abdess@orphee.polytechnique.fr}}
\vskip 1truecm
\centerline{\it Centre de Physique Th{\'e}orique, Ecole Polytechnique}  
\centerline{\it 91128 Palaiseau Cedex, France.}
\centerline{\it Laboratoire Propre du CNRS UPR A.0014.}
\smallskip
\vskip 1truecm


\noindent
{\bf Abstract}

We show how to adapt the Gelfand-Zetlin basis for describing the atypical
representation of ${\cal U}_{\displaystyle{q}}(sl(N))$ when $q$ is root of 
unity. The explicit construction of atypical representation is presented in 
details for $N=3$.

\vskip 1.5truecm

\line{\hfill                                                {\it September 95}}

\eject

{\bfw 1. Introduction.}

The present paper is a sequel to $[1]$, in which we presented an improvement of
$[2]$ giving Gelfand--Zetlin construction of irreducible representation of 
${\cal U}_{\displaystyle{q}}(sl(N))$ at roots of unity independently of their 
nature. We have shown that it is possible to describe the periodic, 
semi-periodic, nilpotent, usual and some atypical representation of 
${\cal U}_{\displaystyle{q}}(sl(N))$ by the fractional parts formalism. 
However, atypical representations generally need a special treatment.

The Gelfand-Zetlin basis in the form $[1]$ is not yet totally adapted 
for atypical representation $[3]$. This has to be compared with the fact that,
for superalgebras, the atypical representations are more difficult to 
describe with the Gelfand--Zetlin than the typical ones: the atypical 
representations of some superalgebras or quantum superalgebras were obtained,
for example, $[4,\;5]$ in the case of $gl(n|1)$ and $[6]$ in the case of
${\cal U}_{\displaystyle{q}}(gl(2|2))$, but the general case is not yet 
treated. 

Note that the paper $[7]$ provides the atypical representations of 
${\cal U}_{\displaystyle{q}}(sl(3))$. The matrix elements of 
${\cal U}_{\displaystyle{q}}(sl(3))$ given by this formalism do not contains 
denominators, and do not generate divergence when $q^{m}=1$. A classification 
of irreducible representations of ${\cal U}_{\displaystyle{q}}(sl(3))$ was 
done in $[13,\;14]$. 

The purpose of this paper is to provide a procedure that enables the 
construction of general atypical representations suitably adapted the 
Gelfand--Zetlin basis.

In section 2, we give a simple example of {\it explicit} construction of the 
{\bf flat} (no multiplicity) representations for ${\cal U}_{\displaystyle{q}}(sl(3))$ and we present the general idea for the construction of atypical 
representation of ${\cal U}_{\displaystyle{q}}(sl(N))$. The atypical case of 
${\cal U}_{\displaystyle{q}}(sl(3))$ is presented in details in section 3.   

{\bfw 2. The Primitive Gelfand--Zetlin basis.}

{\bfw 2.1. The quantum algebra ${\cal U}_{\displaystyle{q}}(sl(N))$.}

The quantum algebra ${\cal U}_{\displaystyle{q}}(sl(N))$ $[8,\;9]$ is defined 
by the generators $k_{i}$, $k_{i}^{-1}$, $e_{i}$, $f_{i}$ ($i=1,\cdots,N-1$) 
and the relations
$$
\eqalign{ 
&k_i e_j = q^{a_{ij}} e_j k_i,\;\;\;\;\;\;\;\; 
k_i f_j  = q^{-a_{ij}} f_j k_i \;,\cr
&[e_i,f_j] = \delta_{ij} {k_i -k_i^{-1} \over q-q^{-1}} \;,\cr
&[e_i,e_j]  = 0 \qquad \hbox{for} \qquad |i-j|>1 \;,\cr
&e_i^2 e_{i\pm 1}  - (q+q^{-1}) e_i e_{i\pm 1} e_i + e_{i\pm 1} e_i^2 = 
0 \;,\cr
&f_i^2 f_{i\pm 1}  - (q+q^{-1}) f_i f_{i\pm 1} f_i + f_{i\pm 1} f_i^2 = 
0 \;,\cr} \eqno(2.1)
$$
The two last equations are called the Serre relations, and 
$( a_{ij})_{i,j=1,...,N-1}$ is the Cartan matrix of $sl(N)$, i.e.
 
$$
a_{ij}=\displaystyle\left\lbrace\matrix{
\hfill 2 &\hbox{for}& \hfill i=j \cr
& & \cr
-1&\hbox{for}& \hfill j=i\pm 1 \cr
& & \cr
0 &\hbox{for}& \hfill |i-j|>1 \cr}\right. \eqno(2.2)
$$

Let us now define the adapted Gelfand--Zetlin basis for the representations 
of ${\cal U}_{\displaystyle{q}}(sl(N))$, the corresponding states are called 
{\bf primitive vectors} of Gelfand-Zetlin pattern.

{\bfw 2.2. Primitive Vectors of the Gelfand--Zetlin basis.}

The states are 
$$
|p\rangle =
\left|
\matrix{p_{1N}& &p_{2N}&\cdots &p_{N-1,N}& &p_{NN} \hfill\cr
&&&&&& \hfill\cr
&p_{1N-1}& &\cdots& &p_{N-1,N-1}& \hfill\cr
&&&&&& \hfill\cr
& & \ddots & \cdots & \adots  & &  \hfill\cr
&&&&&& \hfill\cr
& &p_{12}& &p_{22}& & \hfill\cr
&&&&&& \hfill\cr
& & &p_{11}& & & \hfill\cr} 
\right.
\hbox{\raise -17 ex\hbox{\psfig{figure=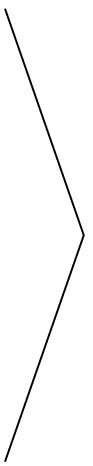}}} \eqno(2.3)
$$
(with respect to $[2]$, we use $p_{ij}=h_{ij}-i$ instead of $h_{ij}$). The 
primitive Gelfand--Zetlin basis $[10]$ is labelled by ${1\over 2}N(N+1)$ 
numbers $p_{ij}$. When $q$ is generic, the first line of indices determines 
the {\it highest weight} of the representation. The states (2.3) within the 
same module $V([p]_{N})$ are distinguished by $p_{ij},\; i,j=1,...,N-1$, which 
assume values consistent with the triangular inequalities
$$
\eqalign{
&p_{i,j+1}-p_{i,j} \in Z_{+}, \cr
&p_{i,j}-p_{i+1,j+1}-1 \in Z_{+},\;\;\;\;\;\;\;\; i,\;j=1,...,N-1, \cr} 
\eqno(2.4)
$$
or,
$$
p_{i,j+1}\geq p_{ij} >  p_{i+1,j+1}. \eqno(2.5)
$$
The dimension of $V([p]_{N})$ is given by
$$
dim V([p]_{N})={\displaystyle\prod_{i=1}^{N-1} \displaystyle\prod_{j=i+1}^{N}
(p_{iN}-p_{jN}) \over\displaystyle\prod_{i=1}^{N-1} (N-i)!}. \eqno(2.6)
$$

\eject

{\bfw 2.3. The primitive representation}

The action of the generators $k_{l}^{\pm 1}$, $e_{l}$ and $f_{l}$ 
$(l=1,2,...,N-1)$ is given by
$$
\eqalign{
&k_{l}^{\pm 1}|p\rangle=q^{\pm (2 \sum_{i=1}^{l} p_{il}-\sum_{i=1}^{l+1}p_{i,l+1}-
\sum_{i=1}^{l-1}p_{i,l-1}-1)} |p\rangle \cr
&f_{l}|p\rangle=\sum_{j=1}^{l}{P_{1}(jl;p)P_{2}(jl;p) \over P_{3}(jl;p)}
|p_{jl}-1\rangle \cr
&e_{l}|p\rangle=\sum_{j=1}^{l}{P_{1}(jl;p_{jl}+1)P_{2}(jl;p_{jl}+1) \over 
P_{3}(jl;p_{jl}+1)}|p_{jl}+1\rangle \cr} \eqno(2.7)
$$
where $|p_{jl}\pm 1 \rangle$ denotes the state differing from $|p\rangle$ by 
only $p_{jl}\rightarrow  p_{jl}\pm 1$, and
$$
\eqalign{
& P_{1}(jl;p)  = \prod_{i=1}^{l+1}
              [\varepsilon_{ij}(p_{i,l+1}-p_{j,l}+1)]^{{1/2}}
,\cr
& P_{2}(jl;p)  = \prod_{i=1}^{l-1}
              [\varepsilon_{ji}(p_{j,l}-p_{i,l-1})]^{{1/2}}
,\cr
& P_3(jl;p)  = \prod_{{i=1 \atop i\ne j}}^{l}
            [\varepsilon_{ij}(p_{i,l}-p_{j,l})]^{1/2}
            [\varepsilon_{ij}(p_{i,l}-p_{j,l}+1)]^{1/2},\cr} \eqno(2.8)
$$
$\varepsilon_{ij}$ being the sign defined by
$$
\varepsilon_{ij}=\left\{\matrix{
\hfill 1&\hbox{for}& \hfill i \leq j \cr
        && \cr 
       -1&\hbox{for}& \hfill i>j \cr
}\right. \eqno(2.9)
$$

In the following, we take $q$ to be a root of unity and $p_{ij}$ are {\it 
integers}. We restrict ourself to the quantum Lie algebra, where the raising 
and lowering operators are nilpotents, i.e. $e_{i}^{m}=f_{i}^{m}=0$ and where 
the Cartan generators $h_{i}$ are such that 
$k_{i}^{m}=\big( q^{h_{i}} \big)^{m}=1$ (representations of this case were 
studied by Lusztig $[15]$). Let $m$ the smallest 
integer such that $q^{m}=1$. We will consider only the case of odd $m$ in 
this paper. A similair discussion is valid when $m$ is even.

We consider here {\it the quantum analogue of classical (highest weight and 
lowest weight) irreducible representations} with a highest weight that obeys
$$
p_{1N}-p_{NN}\;>\;m, \eqno(2.14)
$$
when $q^{m}=1$ these representations are not always irreducible, since some
new {\it singular vectors} arise in the corresponding Verma module, that are 
not obtained from the highest weight vector by action of the translated Weyl 
group. Quotienting by the sub-representation generated by the singular vector 
leads to new irreducible representations that we call {\it atypical by 
analogy with the case of superalgebras}. 

\eject

{\bfw 2.4. A example of flat representation}

In $[1]$, we have introduced some parameters to {\it break the symmetry 
between the actions of $e_{l}$ and $f_{l}$}. They were taken to be $0$, $1/2$ 
or $1$. A good choise of these parameters permits the elimination of the 
singular vectors (these singular vectors are states arising in the r.h.s of
$(2.7)$ but do not obey to the triangular inequalities) in a natural way. If 
an equal number of factors in numerators and denominators are simulaneously 
equal to zero, and if the vector from r.h.s of $(2.7)$ is a singular vector, 
we can ajust these parameters such that the number of zeroes in numerator is 
superior to the number of zeroes in denominator. This procedure describes 
successfully the {\bf flat} representations $[1,\;11,\;12,\;13,\;14]$, i.e. 
when
$$
p_{1N}-p_{NN}=m+1. \eqno(2.10)
$$
For example, the representations of ${\cal U}_{\displaystyle{q}}(sl(3))$ of 
dimension $7$ for $m=3$ ($p_{13}=4$, $p_{23}=2$, $p_{33}=0$) is described by
$$
\eqalign{
&f_{1}|p\rangle = \biggl([p_{11}-p_{22}-1]\biggr)^{1/2}|p_{11}-1\rangle, \cr
&f_{2} |p\rangle =\biggl({[p_{13}-p_{12}+1][p_{12}-p_{23}-1][p_{12}-p_{33}-1]
\over [p_{12}-p_{22}-1][p_{12}-p_{22}]}\biggr)^{1/2}[p_{12}-p_{11}]
|p_{12}-1\rangle \cr
&\;\;\;\;\;\;\;\;+\biggl({[p_{13}-p_{22}+1][p_{23}-p_{22}+1][p_{11}-p_{22}]
\over [p_{12}-p_{22}+1][p_{12}-p_{22}]}\biggr)^{1/2}[p_{22}-p_{33}-1]
|p_{22}-1\rangle, \cr 
&e_{1}|p\rangle = [p_{12}-p_{11}]\biggl([p_{11}-p_{22}]\biggr)^{1/2}
|p_{11}+1\rangle, \cr
&e_{2} |p\rangle =\biggl({[p_{13}-p_{12}][p_{12}-p_{23}][p_{12}-p_{33}]
\over [p_{12}-p_{22}+1][p_{12}-p_{22}]}\biggr)^{1/2}|p_{12}+1\rangle \cr
&\;\;\;\;\;\;\;\;+\biggl({[p_{13}-p_{22}][p_{23}-p_{22}][p_{11}-p_{22}-1]
\over [p_{12}-p_{22}-1][p_{12}-p_{22}]}\biggr)^{1/2}|p_{22}+1\rangle, \cr}
\eqno(2.11)
$$ 
we remark that 
$$
\eqalign{
&f_{1}\left|
\matrix{4& &2& &0 \hfill \cr
&3& &2& \hfill \cr
& &3 & &  \hfill \cr}\right\rangle=f_{2}\left|
\matrix{4& &2& &0 \hfill \cr
&3& &2& \hfill \cr
& &3 & &  \hfill \cr}\right\rangle=0, \cr}\eqno(2.12)
$$
and
$$
\eqalign{
&e_{1}\left|
\matrix{4& &2& &0 \hfill \cr
&3& &2& \hfill \cr
& &3 & &  \hfill \cr}\right\rangle=e_{2}\left|
\matrix{4& &2& &0 \hfill \cr
&3& &2& \hfill \cr
& &3 & &  \hfill \cr}\right\rangle=0, \cr}\eqno(2.13)
$$
the others states form just the irreducible representation of 
dimension $7$ ({\it no multiplicity}). A similair methode describe also the 
representation of dimension $18$ for $m=5$ ($p_{13}=6$, $p_{23}=2$, 
$p_{33}=0$) (see the figure in $[11]$).  

\medskip
\centerline{\psfig{figure=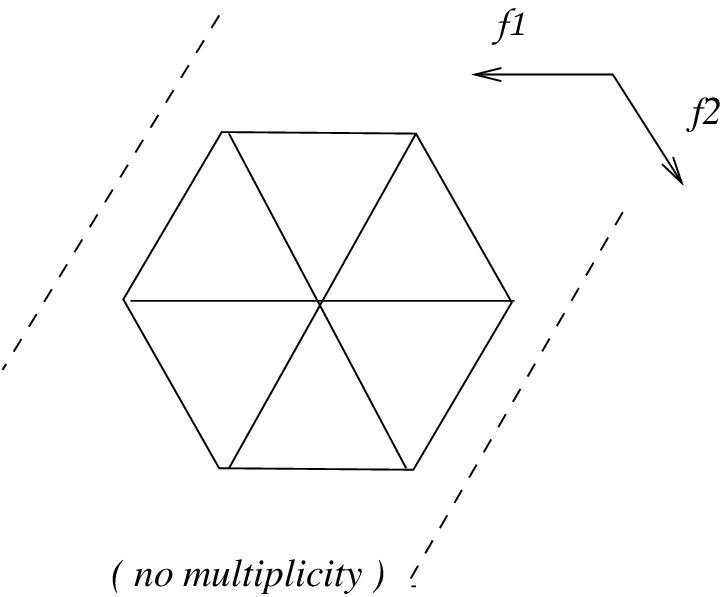}}
\medskip
\centerline{\it Fig. 1: representation of ${\cal U}_{\displaystyle{q}}(sl(3))$
of dimension $7$ for $m=3$. }
\bigskip

{\bfw 2.5. Atypical representations and adaptation of the Gelfand--Zetlin 
basis} 

Let $\eta_{jl}$ and $\eta'_{jl}$ be respectively the numbers of zero of 
$P_{1}(jl;p)\;P_{2}(jl;p)$ and $P_{3}(jl;p)$. We note that the maximum value 
of $\eta'_{jl}$ is $l-1$. If a primitive vector from the r.h.s of $(2.7)$ does
not belong to the module under consideration, then the corresponding term is 
zero ($\eta_{jl} >\eta'_{jl}$) . If an equal numbers of factors in numerators 
and denominators are simultaneously equal to zero ($\eta_{jl} =\eta'_{jl}$), 
they should be cancelled out and the corresponding primitive vector here 
belongs to the module. If $\eta'_{jl} > \eta_{jl}$, the matrix elements of 
$f_{l}$ are undefined. In the next, we will give how to eliminate the 
divergences of the matrix elements using change of bases. This method will be 
illustrated on a simple example in section 3.

Suppose     
$$
\eta_{jl}=0,\;\;\;\;\;\;\hbox{and}\;\;\;\;\;\;\eta'_{jl}=l-1,\;\;\;\;1\leq j 
\leq l  \eqno(2.15)
$$
i.e.
$$
\eqalign{
&p_{i,l+1}-p_{jl}+1\not =0\;[\;m\;],\;\;\;\;\;\;\;\;1\leq i \leq l+1,
\;\;\;\;\;\;\;\;1\leq j 
\leq l, \cr
&p_{jl}-p_{i,l-1}\not =0\;[\;m\;],\;\;\;\;\;\;\;\;1\leq j \leq l,
\;\;\;\;\;\;\;\;1\leq i \leq l-1, \cr 
&p_{il}-p_{jl}=0\;[\;m\;],\;\;\;\;\;\;\;\;1\leq j \leq l,
\;\;\;\;\;\;\;\;1\leq i \leq l, \cr} \eqno(2.16)
$$
thus, there exists $\beta_{1},\;\beta_{2},\cdots , \;\beta_{l-1} \in Z_{+}$ such that
$$
\eqalign{
&p_{i,l}-p_{i+1,l}=\beta_{i}\;m,\;\;\;\;\;\;\;\;1\leq i \leq l-1 \cr
&p_{i,l}-p_{j,l}=(\beta_{i}+\beta_{i+1}+\cdots+\beta_{j-1})\;m,
\;\;\;\;\;\;\;\;\;\;j > i. \cr} 
\eqno(2.17)
$$
The action of $f_{l}$ over a state satisfying $(2.16)$ produces in the r.h.s 
of $(2.7)$ a set of states $\lbrace |p'_{1l}\cdots p'_{il}\cdots p_{ll}\rangle,\;1 \le  i \le l \rbrace$, where 
$$
\eqalign{
&p'_{1l}=\left\{\matrix{
\hfill p_{il}+(\beta_{1}+\cdots+\beta_{i-1})\;m&\hbox{for}& \hfill i > 1 \cr
        && \cr
       p_{1l}-1&\hbox{for}& \hfill i=1 \cr
}\right. \cr
&p'_{il}=\left\{\matrix{
\hfill p_{1l}-1-(\beta_{1}+\cdots+\beta_{i-1})\;m&\hbox{for}& \hfill i > 1 \cr
        && \cr
       p_{1l}-1&\hbox{for}& \hfill i=1 \cr
}\right. \cr} \eqno(2.18)
$$
this set is called a set of {\bf type $p_{1l}-1$} {\bf or} 
$(\lbrace 1\rbrace ,\;\lbrace 2, \cdots,\;l \rbrace)$. 

\proclaim Definition.1. Let a state satisfying 
$$
\eqalign{
&p_{1l}-p_{il}\pm \beta=0\;[\;m\;],\;\;\;\;\;\;\;i=2,\cdots,l,\;\;\;\;\;\;
1 \le |\beta | \le \;m-1, \cr
&p_{il}-p_{jl}=0\;[\;m\;],\;\;\;\;\;\;\;i,\;j=2,\cdots ,\;l. \cr}
\eqno(2.19)
$$
This state is called a state of {\bf type} $p_{1l}$ {\bf or} 
$(\lbrace 1\rbrace ,\;\lbrace 2, \cdots,\;l \rbrace)$. Define the operation 
$$
\pi_{1i}^{l} ( |p_{1l}\cdots p_{il}\cdots p_{ll}\rangle )=
|p'_{1l}\cdots p'_{il}\cdots p_{ll}\rangle   \eqno(2.20)
$$
where
$$
p'_{1l}=\left\{\matrix{
\hfill p_{il}+(\beta_{1}+\cdots+\beta_{i-1})\;m&\hbox{for}& \hfill i > 1 \cr
        && \cr
       p_{1l}&\hbox{for}& \hfill i=1 \cr
}\right. \eqno(2.21)
$$
$$
p'_{il}=\left\{\matrix{
\hfill p_{1l}-(\beta_{1}+\cdots+\beta_{i-1})\;m&\hbox{for}& \hfill i > 1 \cr
        && \cr
       p_{1l}&\hbox{for}& \hfill i=1 \cr
}\right. \eqno(2.22)
$$
This operation is called {\bf exchange mapping of level $l$ centered on 1}, and
the set $\lbrace |p'_{1l}\cdots p'_{il}\cdots p_{ll}\rangle,\;1 \le  i \le l \rbrace$ is of {\bf type} $p_{1l}$ {\bf or} 
$(\lbrace 1\rbrace ,\;\lbrace 2, \cdots,\;l \rbrace)$. This set of states 
has the same eigenvalues for the Cartan operators (degenerate states) and has
to satisfy the triangular inequalities.

\proclaim Lemma. Let a state satisfy the condition $(2.19)$, and let 
$\lbrace |p'_{1l}\cdots p'_{il}\cdots p_{ll}\rangle,\;1 \le  i \le l \rbrace$
be the set of all states obtained by action of the mapping $\pi_{1i}^{l}$ 
over this state. This set is {\bf isomorphic} to a set obtained by action of
the mapping $\pi_{\mu \nu}^{l}$ over a state of type {\bf type 
$p_{\mu l}\mp \beta$}$\;\;(\mu \not = 1)$. 

\proclaim Proof. Let a state be of type $(\lbrace 1 \rbrace , \lbrace 2, \cdots, \;l\rbrace$, 
i.e.
$$
\eqalign{
&p_{1l}-p_{il}\pm \beta= 0\;[\;m\;]\;\;\;\;\;\;\;i=2,\cdots,l,\;\;\;\;\;\;
1 \le |\beta | \le \;m-1, \cr
&p_{il}-p_{jl}=0\;[\;m\;]\;\;\;\;\;\;\;i\not=1\;\;\;\;\hbox{and}
\;\;\;\;j\not=1, \cr}
$$
and, let
$$
\pi_{1\mu}^{l} ( |p_{1l}\cdots p_{\mu l}\cdots p_{ll}\rangle )=
|p_{1l}\pm \beta\cdots p_{\mu l}\mp \beta \cdots p_{ll}\rangle   
$$
The state in the r.h.s is of {\bf type $p_{\mu l}\mp \beta$}$\;\;(\mu \not = 1)$, and
$$
\pi_{\mu \nu}^{l} \circ \pi_{1 \mu}^{l}=\pi_{1 \nu}^{l}. \eqno(2.23)
$$  

All states of the set $(2.18)$ are obtained by action of the mapping 
$\pi_{1i}^{l}$ over the state $|p_{1l}-1\cdots p_{il}\cdots p_{ll}\rangle$.
Using this notation, the action of $f_{l}$ is  
$$
\eqalign{
&f_{l}|p\rangle = \sum_{j=1}^{l} {P_{1}(jl;\;p)\;P_{2}(jl;\;p) \over 
P_{3}(jl;\;p)}|p'_{1l}\cdots p'_{jl}\cdots p_{ll}\rangle \cr}
\eqno(2.24)
$$
We note that the number of zero in the polynomials $P_{3}(jl;\;p)$ is $l-1$.

\proclaim Definition.2. Let a set be of state of {\bf type} $(\lbrace 1\rbrace ,\;\lbrace 2, \cdots,\;l \rbrace)$, i.e. 
$$
p_{1l}-p_{il}\pm \beta =0\;[\;m\;],\;\;\;\;\;\;\;\;\;\;2\le i \le l,\;\;\;\;\;\;
1 \leq |\beta \leq m-1,
$$
and, let the new basis be given by
$$
|p'_{1l}...p'_{il}...p_{ll}\rangle=\sum_{j=1}^{l}\;D_{ij}\;
\parallel p'_{1l}...p'_{jl}...p_{ll}\rangle \eqno(2.25)
$$
where, $D$ is a $l \times l$ rotation matrix, i.e.
$$
D^{t}. D=D. D^{t}=\hbox{{\bf 1}}. \eqno(2.26)
$$
This new basis is called the {\bf modified basis of type 
$(\lbrace 1\rbrace ,\;\lbrace 2, \cdots,\;l \rbrace)$}. The 
primitive and the modified sets satisfy the triangular inequalities $(2.5)$. 

The finiteness of the matrix elements $\langle p|e_{l}f_{l}|p\rangle$ and 
$\langle p|f_{l}e_{l}|p\rangle$ (preserve 
$[e_{l},\;f_{l}]={k_{l}-k_{l}^{-1} \over q-q^{-1}}$) imply that there exists 
a modified basis such that the new matrix elements are without 
divergences. Using this definition, the equation $(2.24)$ is reduced to
$$
\eqalign{
f_{l}|p\rangle &= \sum_{i=1}^{l}\biggl(\sum_{j=1}^{l} 
{P_{1}(jl;\;p)\;P_{2}(jl;\;p) \over P_{3}(jl;\;p)} D_{ji}\biggr)\parallel p'_{1l}\cdots p'_{il}
\cdots p_{ll}\rangle \cr
&=\sum_{i=1}^{l} A_{il} \parallel p'_{1l}\cdots p'_{il}
\cdots p_{ll}\rangle \cr}
\eqno(2.27)
$$
where, $A_{il}$ are the new matrix elements associated to the modified 
basis, i.e. 
$$
\eqalign{
&A_{il}=\sum_{j=1}^{l} 
{P_{1}(jl;\;p)\;P_{2}(jl;\;p) \over P_{3}(jl;\;p)} D_{ji},\;\;\;\;\;\;\;1\le i \le l, \cr}
\eqno(2.28)
$$
this matrix elements are called the {\bf modified matrix elements}. Generally, 
we choose  
$$
\eqalign{
&A_{1l}=\biggl(\sum_{j=1}^{l} {P^{2}_{1}(jl;\;p)\;P^{2}_{2}(jl;\;p) \over P^{2}_{3}(jl;\;p)} \biggr)^{1/2}, 
\cr
&A_{il}=0\;\;\;\;\;\;\;\;\;\;\;2 \le i \le l, \cr} \eqno(2.29)
$$ 
i.e.
$$
\eqalign{
&f_{l}|p\rangle =\biggl(\sum_{j=1}^{l} {P^{2}_{1}(jl;\;p)\;P^{2}_{2}(jl;\;p) \over P^{2}_{3}(jl;\;p)} \biggr)^{1/2} \parallel p_{1l}-1\cdots p_{ll}\rangle . 
\cr} \eqno(2.30)
$$
We note that the matrix element of $(2.30)$ is finit (i.e. without 
divergences).

We are now able to claim this generalization,

\proclaim Definition.3. Let a collection $\lbrace \hbox{\bf I}_{k},\; k\in 
\hbox{\bf J}\subset \hbox{\bf N}\rbrace$ of subsets of 
$\lbrace 1,\cdots,\;l\rbrace$ satisfying the following conditions:
$$
\eqalign{
&\hbox{\bf 1.}\;\;\;\;\;\;\;\;\;\;\;\;\;\;\;\;\;\;\;\;\;\;\;\;\;\;\;\;\;\;\;\;\;\;\hbox{ if}\;\;\;k < s\;\;\;\;\;\;\forall i \in\hbox{\bf I}_{k}\;\;\; \forall j \in\hbox{\bf I}_{s}\;\;\;\hbox{thus}\;\;\;\; i<j, \cr
&\hbox{\bf 2.}\;\;\;\;\;\;\;\;\;\;\;\;\;\;\;\;\;\;\;\;\;\;\;\;\;\;\;\;\;\;\;\;\;\;\hbox{\bf I}_{k} \;\cap\; \hbox{\bf I}_{s} = \emptyset,\;\;\;\;\;\;\;\hbox{if}\;\;\;\;\;\;\;k\not = s,  \cr
&\hbox{\bf 3.}\;\;\;\;\;\;\;\;\;\;\;\;\;\;\;\;\;\;\;\;\;\;\;\;\;\;\;\;\;\;\;\;\;\;\displaystyle\bigcup_{k\in\;\hbox{\bf J}}\;\;\hbox{\bf I}_{k}=\lbrace 1,\cdots,\;l\rbrace, \cr }
\eqno(2.31)
$$
and, where
$$
\eqalign{
&p_{il}-p_{jl}=\;[\;m\;],\;\;\;\;\;\;i,j \in \hbox{\bf I}_{k},\;\;\;\;\;i< j, 
\cr
&p_{il}-p_{jl} = \zeta_{ks}\;[\;m\;]\;\;\;\;\;\;i \in \hbox{\bf I}_{k}, 
\;j \in 
\hbox{\bf I}_{s}\;\;\;\;\;\;\;k < s,\;\;\;\;\;1 \le |\zeta_{ks} | \le \;m-1. 
\cr} \eqno(2.32)
$$
This state is called a state of {\bf type} $(\hbox{\bf I}_{k},\;k\in\hbox{\bf 
J}\subset \hbox{\bf N})$. Let $\pi^{l}_{ks;ij}$ be the {\bf exchange 
mapping of level $l$ between the subsets $\hbox{\bf I}_{k}$ and 
$\hbox{\bf I}_{s}$} ($k<s$), i.e. if
$$
p_{il}-p_{jl}=\zeta_{ks}+(\beta_{i}+\beta_{i+1}+\cdots+\beta_{j-1})\;m
\eqno(2.33)
$$
we have
$$
\pi^{l}_{ks;ij} ( |p_{1l}\cdots p_{il}\cdots p_{jl}\cdots p_{ll}\rangle )=
|p_{1l}\cdots p'_{il}\cdots p'_{jl}\cdots p_{ll}\rangle   
\eqno(2.34)
$$
with
$$
\eqalign{
&p'_{il}=p_{jl}+(\beta_{i}+\beta_{i+1}+\cdots+\beta_{j-1})\;m=p_{il}+\zeta_{ks}, \cr
&p'_{jl}=p_{il}-(\beta_{i}+\beta_{i+1}+\cdots+\beta_{j-1})\;m=p_{jl}-\zeta_{ks}, \cr}
\eqno(2.35)
$$
The operators $k_{i}$ commutes with exchange mapping $\pi^{l}_{ks;ij}$. Let
$|p\rangle$ be state satisfying $(2.32)$ and $(2.33)$. Define
$$
V_{l}(\hbox{\bf I}_{k},\;k\in \hbox{\bf J})=
\big\lbrace \prod_{\scriptstyle k_{\alpha}\not = s_{\alpha} 
\atop \scriptstyle i_{\alpha} \not = j_{\alpha}}
\pi_{k_{\alpha}s_{\alpha};i_{\alpha}j_{\alpha}}
( |p_{1l}\cdots p_{il}\cdots p_{jl}\cdots p_{ll}\rangle )\;\;\;\;\big\rbrace
\eqno(2.36)
$$
The set $V_{l}\;(\hbox{\bf I}_{k},\;k\in \hbox{\bf J}\subset \hbox{\bf N})$ 
is obtained by the all possible changes between the different subsets 
$\hbox{\bf I}_{k}\; (k\in \hbox{\bf J})$. We note that the all states of 
$V_{l}\;(\hbox{\bf I}_{k},\;k\in \hbox{\bf J}\subset \hbox{\bf N})$
have the same eigenvalues for the Cartan operators.

For example,
$$
\eqalign{
&\hbox{J}=\lbrace 1 \rbrace,\;\;\;\;\;\;\;\;\;\;
\hbox{\bf I}_{1}=\lbrace 1,\cdots ,l \rbrace,\;\;\;\;\;\;\;\;\;\;\;\;\;\;\;\;\;\;\;\;\;\;\;\;\;\;\;\;\;\;\;\;\;\;\;\;\;\;\;\;\;\;\;\;dim V_{l}=1, \cr
&\hbox{J}=\lbrace 1,\;2 \rbrace,\;\;\;\;\;\;\hbox{\bf I}_{1}=\lbrace 1 \rbrace,\;\;\;\;\;\;\;\;\;\;\;\;\;\;\;\;\hbox{\bf I}_{2}=\lbrace 2,\cdots ,\;l \rbrace,\;\;\;\;\;\;\;\;\;\;\;\;dim V_{l} = l, \cr
&\hbox{J}=\lbrace 1,\;2 \rbrace,\;\;\;\;\;\;\hbox{\bf I}_{1}=\lbrace 1,\;2 \rbrace,\;\;\;\;\;\;\;\;\;\;\;\;\hbox{\bf I}_{2}=\lbrace 3,\cdots ,\;l \rbrace,\;\;\;\;\;\;\;\;\;\;\;\;dim V_{l}= {1\over 2}l(l-1). \cr}
$$

\proclaim Definition.4. Let a new basis given by 
$$
|p'_{1l}\cdots p'_{il}\cdots p_{ll}\rangle=\sum_{j=1}^{l}\;D_{ij}\;
\parallel p'_{1l}\cdots p'_{jl}\cdots p_{ll}\rangle \eqno(2.37)
$$
where the states $|p'_{1l}\cdots p'_{il}\cdots p_{ll}\rangle$ are in the set 
$V_{l}\;(\hbox{\bf I}_{k},\;k\in \hbox{\bf J}\subset \hbox{\bf N})$ and
$D$ is a $dim\;V_{l} \times dim\;V_{l}$ rotation matrix. This new basis is 
called the {\bf modified basis of type}
$(\hbox{\bf I}_{k},\;k\in\hbox{\bf J}\subset \hbox{\bf N})$ . We note that 
the primitive and the modified bases satisfy the triangular inequalities. The 
set of states $\parallel p'_{1l}\cdots p'_{jl}\cdots p_{ll}\rangle$ is noted
by $\overline{V}_{l}\;(\hbox{\bf I}_{k},\;k\in \hbox{\bf J}\subset 
\hbox{\bf N})$.  
$$
\diagram{
(\lbrace 1,\cdots ,\; l \rbrace ) & \hfl{f_{l}}{} & (\lbrace 1\rbrace,\lbrace 2, \;l \rbrace ) 
& \hfl{f_{l}}{} & (\lbrace 1\rbrace,\lbrace 2,\;l \rbrace ) & \cdots  \cr
&& \vfl{f_{l}}{}&& \vdots & \cr
&& (\lbrace 1,\;2 \rbrace, \lbrace 3,\cdots,\;l \rbrace ) &\cdots & &\cr 
&& \vfl{f_{l}}{} & & &\cr 
&& (\lbrace 1,\;2,\;3 \rbrace, \lbrace 4,\cdots,\;l \rbrace ) &\cdots & &\cr 
&& \vdots & & &\cr }
$$
\centerline{\it Fig. 2: different transition between the type of sets in the 
case of ${\cal U}_{\displaystyle{q}}(sl(N))$. }

Using the definition {\bf 3}, the action of $f_{l}$ over a modified state of 
{\bf type} ($\lbrace 1\rbrace$, $\lbrace 2, \cdots,\;l \rbrace$) produces  
two sets of states respectively of {\bf type} ($\lbrace 1\rbrace$, 
$\lbrace 2, \cdots,\;l \rbrace$) and ($\lbrace 1,\;2\rbrace$, 
$\lbrace 3, \cdots,\;l \rbrace$). The matrix elements of 
($\lbrace 1\rbrace$, $\lbrace 2, \cdots,\;l \rbrace$) $\rightarrow$ 
($\lbrace 1\rbrace$, $\lbrace 2, \cdots,\;l \rbrace$) do not 
contain divergences (see the structure of $P_{3}(jl;\;p)$). But if the 
matrix elements of ($\lbrace 1\rbrace$, $\lbrace 2, \cdots,\;l \rbrace$)
$\rightarrow$ ($\lbrace 1,\;2\rbrace$, $\lbrace 3, \cdots,\;l \rbrace$) 
contain somes divergences, using of the definition {\bf 4} we take a rotation 
in this second set in such a way as to eliminate these divergences 
( $\langle p|e_{l}f_{l}|p\rangle$ and $\langle p|f_{l}e_{l}|p\rangle$ have to 
remain finit). We have to repeat this mechanism as far as the elimination 
of the all divergences of the Gelfand--Zetlin representation.  

{\bfw 3. Applications.}
 
In the following, we present in details the explicit construction of atypical 
representation for $N=3$. We will consider only the operator
$f_{1}$ and $f_{2}$, a similar discussion is valid for $e_{1}$ and 
$e_{2}$. The primitive Gelfand--Zetlin state for this particular case is just
$$
|p\rangle =
\left|
\matrix{p_{13}& &p_{23}& &p_{33} \hfill \cr
&p_{12}& &p_{22}& \hfill \cr
& &p_{11} & &  \hfill \cr} \right\rangle \eqno(3.1) 
$$
where $p_{33}$ is choosen equal to zero. The actions of the generators 
$f_{1}$ and $f_{2}$ for ${\cal U}_{\displaystyle{q}}(sl(3))$ are given by
$$
f_{1} |p\rangle = \biggl([p_{12}-p_{11}+1][p_{11}-p_{22}-1]\biggr)^{1/2} 
|p_{11}-1\rangle \eqno(3.2) 
$$
$$
\eqalign{
&f_{2} |p\rangle =\biggl({[p_{13}-p_{12}+1][p_{12}-p_{23}-1][p_{12}-p_{33}-1]
[p_{12}-p_{11}] \over [p_{12}-p_{22}][p_{12}-p_{22}-1]}\biggr)^{1/2}
|p_{12}-1\rangle , \cr
&\;\;\;\;\;\;\;\;\;+\biggl({[p_{13}-p_{22}+1][p_{23}-p_{22}+1][p_{22}-p_{33}-1]
[p_{11}-p_{22}] \over [p_{12}-p_{22}][p_{12}-p_{22}+1]}\biggr)^{1/2}
|p_{22}-1\rangle. \cr} \eqno(3.3)
$$ 
We note that in this case there are only two type of sets 
$(\lbrace 1,\;2 \rbrace )$ and $(\lbrace 1\rbrace,\;\lbrace 2\rbrace)$ and 
the maximal number of zero in the denominator of the matrix
elements is only one.

Now suppose that  
$$
\eqalign{
&p_{i3}-p_{j2}+1 \not\equiv 0\;[\;m\;],\;\;\;\;1\leq i \leq 3,\;\;\; 
\hbox{and} \;\;\;1\leq j \leq 2, \cr
&p_{i2}-p_{11} \not\equiv 0\;[\;m\;],\;\;\;\;1\leq i \leq 2, \cr} \eqno(3.4)
$$
i.e. 
$$
\eta_{12}=\eta_{22}=0.  \eqno(3.5)
$$
The matrix elements of $f_{2}$ are infinite if 
$$
\eqalign{
&(a)\;\;\;\;\;\;\;\;\;\;\;\;p_{12}-p_{22}=0\;[\;m\;], \cr
&(b)\;\;\;\;\;\;\;\;\;\;\;\;p_{12}-p_{22}+1=0\;[\;m\;], \cr
&(c)\;\;\;\;\;\;\;\;\;\;\;\;p_{12}-p_{22}-1=0\;[\;m\;]. \cr} \eqno(3.6)
$$
$$
\diagram{
(\lbrace 1,\;2 \rbrace ) & \hfl{f_{2}}{} & (\lbrace 1\rbrace,\lbrace 2 \rbrace ) 
& \hfl{f_{2}}{} & (\lbrace 1\rbrace,\lbrace 2 \rbrace ) & \cdots  \cr
&& \vfl{f_{2}}{}&& \vdots & \cr
&& (\lbrace 1,\;2 \rbrace ) & & &\cr 
&& \vfl{f_{2}}{} & & &\cr 
&&  (\lbrace 1\rbrace,\lbrace 2 \rbrace )& \cdots & &\cr
&&  \vdots & & &\cr}
$$
\centerline{\it Fig. 3: different transition between the type of sets 
$(\lbrace 1,\;2 \rbrace )$ and }
\centerline{\it $(\lbrace 1\rbrace,\lbrace 2 \rbrace )$  in the 
case of ${\cal U}_{\displaystyle{q}}(sl(3))$. }

\proclaim Case (a). Let a state satisfy the following conditions
$$
p_{12}-p_{22}=0\;\;[\;m\;]\;\;\;\;\;\hbox{or}\;\;\;\;\;\;\;p_{12}-p_{22}=\beta\;\;m\;\;\;\;\;\;\beta \in Z_{+} \eqno(3.7)
$$
i.e. 
$$
\eta'_{12}=\eta'_{22}=1.  \eqno(3.8)
$$
The action action of $f_{2}$ on this state give 
$$
f_{2}|p_{12}\;\;p_{22}\rangle= {\kappa \over([\;m\;][m-1])^{1/2}} |p_{12}-1\;\;
p_{22}\rangle + {\kappa \over ([\;m\;][m+1])^{1/2}} |p_{12}\;\;p_{22}-1\rangle ,
\eqno(3.9)
$$
where 
$$
\kappa=\biggl([p_{13}-p_{22}+1][p_{23}-p_{22}+1][p_{22}-p_{33}-1][p_{11}-
p_{22}]\biggr)^{1/2}. \eqno(3.10)
$$

Using the definitions {\bf 3} and {\bf 4}, the relation between the primitive 
and the modified states is 
$$
\pmatrix{|p_{12}-1\;\;p_{22}\rangle \cr
\cr
|p_{22}+\beta\;m\;\;p_{12}-1-\beta\;m\rangle \cr}=
D(\phi)\;\pmatrix{\parallel p_{12}-1\;\;p_{22}\rangle \cr
\cr
\parallel p_{22}+\beta\;m\;\;p_{12}-1-\beta\;m\rangle \cr}, 
$$
i.e.
$$
\pmatrix{|p_{12}-1\;\;p_{22}\rangle \cr
\cr
|p_{12}\;\;p_{22}-1\rangle \cr}=
D(\phi)\;
\pmatrix{\parallel p_{12}-1\;\;p_{22}\rangle \cr
\cr
\parallel p_{12}\;\;p_{22}-1\rangle \cr}, \eqno(3.11)
$$
for any values of $p_{11}$ satisfying the triangular inequality respectively 
for the primitive and the modified basis, and where
$$
D(\phi)=\pmatrix{\cos \phi& & \sin \phi \cr
& & \cr
-\sin \phi& & \cos \phi \cr} \eqno(3.12)
$$
Finally,
$$
\eqalign{
f_{2}|p_{12}\;\;p_{22}\rangle &=\kappa \biggl({ \cos\phi \over ([\;m\;][m-1])^{1/2}}-
{\sin\phi \over ([\;m\;][m+1])^{1/2}}\biggr)\parallel p_{12}-1\;\;p_{22}\rangle \cr
&+\kappa \biggl({ \sin\phi \over ([\;m\;][m-1])^{1/2}}+
{\cos\phi \over ([\;m\;][m+1])^{1/2}}\biggr)\parallel p_{12}\;\;p_{22}-1\rangle \cr} \eqno(3.13)
$$

Using the trivial identity,
$$
\eqalign{
&{[a+1] \over [2][a]}+{[a-1] \over [2][a]} =1 \cr} \eqno(3.14)
$$
choose 
$$
\eqalign{
&\cos\phi=\biggl({[m-1] \over [2][\;m\;]}\biggr)^{1/2},\;\;\;\;\;\;\;\;
\sin\phi=\biggl({[m+1] \over [2][\;m\;]}\biggr)^{1/2}, \cr} \eqno(3.15)
$$
The r.h.s of (3.13) is reduced to 
$$
f_{2}|p_{12}\;\;p_{22}\rangle = 
\biggl({[2][p_{13}-p_{22}+1][p_{23}-p_{22}+1][p_{22}-p_{33}-1][p_{11}-
p_{22}] \over [p_{12}-p_{22}-1][p_{12}-p_{22}+1]}\biggr)^{1/2} 
\parallel p_{11}\;\;p_{22}-1\rangle . \eqno(3.16)
$$
This equation correspond to the transition $(\lbrace 1,\;2 \rbrace)\rightarrow
(\lbrace 1\rbrace ,\;\lbrace 2 \rbrace)$.

In this case the action of $f_{1}$ over the modified Gelfand--Zetlin basis is 
given by
$$
f_{1}\parallel p_{12}-1\rangle=\left\{\matrix{
\biggl([p_{12}-p_{11}][p_{11}-p_{22}-1]\biggr)^{1/2}\parallel p_{12}-1\rangle_{p_{11}-1}&
\hbox{if}& p_{11}\not=p_{22}+1 \cr
& & \cr
& & \cr
[m-1]\biggl({[m+1]\over [2]}\biggr)^{1/2}|p_{22}-1\rangle_{p_{11}-1}&
\hbox{if}& p_{11}=p_{22}+1 \cr}\right. \eqno(3.17)
$$
and 
$$
f_{1}\parallel p_{22}-1\rangle=\left\{\matrix{
\biggl([p_{12}-p_{11}+1][p_{11}-p_{22}]\biggr)^{1/2}\parallel p_{22}-1\rangle_{p_{11}-1}&
\hbox{if}& p_{11}\not=p_{22}+1 \cr
& & \cr
& & \cr
\biggl({[m-1]\over [2]}\biggr)^{1/2}|p_{22}-1\rangle_{p_{11}-1}&
\hbox{if}& p_{11}=p_{22}+1 \cr}\right. \eqno(3.18)
$$

\proclaim Remark. If
$$
p_{12}-p_{22}+1=0\;[\;m\;]\;\;\;\;\;\;\hbox{i.e.}
\;\;\;\;\;\;p_{12}-p_{22}+1=\beta\;m\;\;\;(\beta\in Z_{+})
\eqno(3.19)
$$
the relation between the primitive and the modified basis is given by
$$
\pmatrix{|p_{12}\;\;p_{22}\rangle \cr
\cr
|p_{12}+1\;\;p_{22}-1\rangle \cr}=
\pmatrix{\cos \phi& & \sin \phi \cr
& & \cr
-\sin \phi& & \cos \phi \cr}
\pmatrix{\parallel p_{12}\;\;p_{22}\rangle \cr
\cr
\parallel p_{12}+1\;\;p_{22}-1\rangle \cr} \eqno(3.20)
$$
where 
$$
\cos\phi =\biggl({[p_{12}-p_{22}] \over [2][p_{12}-p_{22}+1]}\biggr)^{1/2},
\;\;\;\;\;\;\;\;\;
\sin\phi =\biggl({[p_{12}-p_{22}+2] \over [2][p_{12}-p_{22}+1]}\biggr)^{1/2}.
\eqno(3.21)
$$

\proclaim Case (b). For a state satisfying the condition
$$
p_{12}-p_{22}+1=0\;[\;m\;]\;\;\;\;\hbox{i.e.}\;\;\;\;\;\;p_{12}-p_{22}+1=\beta\;m\;\;
(\beta\in Z_{+}) \eqno(3.22)
$$
the correspondance between the primitive vectors and the modified vectors is 
given by the following formulas
$$
\pmatrix{\parallel p_{12}\;\;p_{22}\rangle \cr
\cr
\parallel p_{12}+1\;\;p_{22}-1\rangle \cr}=
\pmatrix{\cos \phi& & -\sin \phi \cr
& & \cr
\sin \phi& & \cos \phi \cr}
\pmatrix{| p_{12}\;\;p_{22}\rangle \cr
\cr
| p_{12}+1\;\;p_{22}-1\rangle \cr} \eqno(3.23)
$$
where $\phi$ is given by (3.21). If we take a extension of the definition of 
the modified formula
$$
\pmatrix{\parallel p_{12}-1\;\;p_{22}\rangle \cr
\cr
\parallel p_{22}+\beta\;m\;\;p_{12}-1-\beta\;m\rangle \cr}=
D(-\phi)\;
\pmatrix{| p_{12}-1\;\;p_{22}\rangle \cr
\cr
| p_{22}+\beta\;m\;\;p_{12}-1-\beta\;m\rangle \cr} 
$$
i.e.
$$
\pmatrix{\parallel p_{12}-1\;\;p_{22}\rangle \cr
\cr
\parallel p_{12}+1\;\;p_{22}-2\rangle \cr}=
\pmatrix{\cos \phi& & -\sin \phi \cr
& & \cr
\sin \phi& & \cos \phi \cr}
\pmatrix{| p_{12}-1\;\;p_{22}\rangle \cr
\cr
| p_{12}+1\;\;p_{22}-2\rangle \cr} \eqno(3.24)
$$
the action of $f_{2}$ over the modified states is reduced to
$$
\eqalign{
&f_{2} \parallel p_{12}\;\;p_{22}\rangle=\cr
&\biggl({[p_{13}-p_{12}+1][p_{12}-p_{23}-1][p_{12}-p_{33}-1]
[p_{12}-p_{11}] \over [p_{12}-p_{22}][p_{12}-p_{22}-1]}\biggr)^{1/2}
\parallel p_{12}-1\;\;p_{22}\rangle \cr} \eqno(3.25)
$$
and
$$
\eqalign{
&f_{2} \parallel p_{12}+1\;\;p_{22}-1\rangle=\cr
&\biggl({[p_{13}-p_{22}+2][p_{23}-p_{22}+2][p_{22}-p_{33}-2]
[p_{11}-p_{22}+1] \over [p_{12}-p_{22}+2][p_{12}-p_{22}+3]}\biggr)^{1/2}
\parallel p_{12}+1\;\;p_{22}-2\rangle+\cr 
&\biggl({[2][p_{13}-p_{22}+1][p_{23}-p_{22}+1][p_{22}-p_{33}-1]
[p_{11}-p_{22}] \over [p_{12}-p_{22}][p_{12}-p_{22}+2]}\biggr)^{1/2}
| p_{12}\;\;p_{22}-1\rangle\cr} \eqno(3.26)
$$
We note that the state $|p_{12}\;\;p_{22}-1 \rangle$ is of type 
($\lbrace 1,\;2 \rbrace$). The equations $(3.25)$ and $(3.26)$ correspond
respectively to the transitions ($\lbrace 1 \rbrace$, $\lbrace 2 \rbrace)
\rightarrow$ ($\lbrace 1 \rbrace,\lbrace 2 \rbrace)$ and 
$(\lbrace 1 \rbrace,\lbrace 2 \rbrace)\rightarrow 
(\lbrace 1 \rbrace,\lbrace 2 \rbrace)+ (\lbrace 1,\; 2 \rbrace)$.   

\proclaim Remark. For example, the action of $f_{2}$ over the extension 
$(3.24)$ gives 
$$
\eqalign{
&f_{2} \parallel p_{12}-1\;\;p_{22}\rangle=\cr
&\biggl({[p_{13}-p_{12}+2][p_{12}-p_{23}-2][p_{12}-p_{33}-2]
[p_{12}-p_{11}-1] \over [p_{12}-p_{22}-1][p_{12}-p_{22}-2]}\biggr)^{1/2}
\parallel p_{12}-2\;\;p_{22}\rangle+\cr 
&\biggl({[p_{13}-p_{22}+1][p_{23}-p_{22}+1][p_{22}-p_{33}-1]
[p_{11}-p_{22}] \over [p_{12}-p_{22}-1][p_{12}-p_{22}]}\biggr)^{1/2}
\parallel p_{12}-1\;\;p_{22}-1\rangle\cr} \eqno(3.27)
$$
where
$$
\pmatrix{\parallel p_{12}-2\;\;p_{22}\rangle \cr
\cr
\parallel p_{22}+\beta\;m\;\;p_{12}-2-\beta\;m\rangle \cr}=
D(-\phi)\;
\pmatrix{| p_{12}-2\;\;p_{22}\rangle \cr
\cr
| p_{22}+\beta\;m\;\;p_{12}-2-\beta\;m\rangle \cr} \eqno(3.28)
$$
$$
\pmatrix{\parallel p_{12}-1\;\;p_{22}-1\rangle \cr
\cr
\parallel p_{22}-1+\beta\;m\;\;p_{12}-1-\beta\;m\rangle \cr}=
D(-\phi)\;
\pmatrix{| p_{12}-1\;\;p_{22}-1\rangle \cr
\cr
| p_{22}-1+\beta\;m\;\;p_{12}-1-\beta\;m\rangle \cr} \eqno(3.29)
$$

\proclaim Case (c). The discussion is similar to the case (b). Let a primitive 
state satisfy the following condition
$$
p_{12}-p_{22}-1=0\;[\;m\;]\;\;\;\;\;\;\hbox{in}\;\;\;\;\;\;p_{12}-p_{22}-1=\beta\;m\;\;\;\;(\beta\in Z_{+}) \eqno(3.30)
$$

Define
$$
\pmatrix{\parallel p_{12}\;\;p_{22}\rangle \cr
\cr
\parallel p_{22}+\beta\;m\;\;p_{12}-\beta\;m\rangle \cr}=
D(-\phi)\;
\pmatrix{| p_{12}-1\;\;p_{22}\rangle \cr
\cr
| p_{22}+\beta\;m\;\;p_{12}-\beta\;m\rangle \cr} 
$$
i.e.
$$
\pmatrix{\parallel p_{12}\;\;p_{22}\rangle \cr
\cr
\parallel p_{12}-1\;\;p_{22}+1\rangle \cr}=
\pmatrix{\cos \phi& & -\sin \phi \cr
& & \cr
\sin \phi& & \cos \phi \cr}
\pmatrix{| p_{12}\;\;p_{22}\rangle \cr
\cr
| p_{12}-1\;\;p_{22}+1\rangle \cr} \eqno(3.31)
$$
and take the extension
$$
\pmatrix{\parallel p_{12}\;\;p_{22}-1\rangle \cr
\cr
\parallel p_{22}-1+\beta\;m\;\;p_{12}-\beta\;m\rangle \cr}=
D(-\phi)\;
\pmatrix{| p_{12}-1\;\;p_{22}-1\rangle \cr
\cr
| p_{22}-1+\beta\;m\;\;p_{12}-\beta\;m\rangle \cr} 
$$
i.e.
$$
\pmatrix{\parallel p_{12}\;\;p_{22}-1\rangle \cr
\cr
\parallel p_{12}-2\;\;p_{22}+1\rangle \cr}=
\pmatrix{\cos \phi& & -\sin \phi \cr
& & \cr
\sin \phi& & \cos \phi \cr}
\pmatrix{| p_{12}\;\;p_{22}-1\rangle \cr
\cr
| p_{12}-2\;\;p_{22}+1\rangle \cr} \eqno(3.32)
$$
where
$$
\cos \phi=\biggl({[p_{12}-p_{22}] \over [2][p_{12}-p_{22}-1]}\biggr)^{1/2},
\;\;\;\;\;\;\;\;\;\;\;\;\sin \phi=
\biggl({[p_{12}-p_{22}-2] \over [2][p_{12}-p_{22}-1]}\biggr)^{1/2}. \eqno(3.33)
$$

We obtain
$$
\eqalign{
&f_{2} \parallel p_{12}\;\;p_{22}\rangle = \cr
&\biggl({[p_{13}-p_{12}+1][p_{12}-p_{23}-1][p_{12}-p_{33}-1]
[p_{12}-p_{11}] \over [p_{12}-p_{22}][p_{12}-p_{22}+1]}\biggr)^{1/2}
\parallel p_{12}\;\;p_{22}-1\rangle \cr} \eqno(3.34)
$$
and
$$
\eqalign{
&f_{2} \parallel p_{12}-1\;\;p_{22}+1\rangle=\cr
&\biggl({[p_{13}-p_{22}+1][p_{23}-p_{22}+1][p_{22}-p_{33}-1]
[p_{11}-p_{22}] \over [p_{12}-p_{22}][p_{12}-p_{22}+1]}\biggr)^{1/2}
\parallel p_{12}-2\;\;p_{22}+1\rangle+\cr
&\biggl({[2][p_{13}-p_{12}+1][p_{12}-p_{23}-1][p_{12}-p_{33}-1]
[p_{12}-p_{11}] \over [p_{12}-p_{22}][p_{12}-p_{22}-2]}\biggr)^{1/2}
| p_{12}-1\;\;p_{22}\rangle\cr} \eqno(3.35)
$$
We note that the state $|p_{12}-1\;\;p_{22} \rangle$ is of type 
($\lbrace 1,\;2 \rbrace$). The equations $(3.34)$ and $(3.35)$ correspond
respectively to the transitions ($\lbrace 1 \rbrace$, $\lbrace 2 \rbrace)
\rightarrow$ ($\lbrace 1 \rbrace,\lbrace 2 \rbrace)$ and 
$(\lbrace 1 \rbrace,\lbrace 2 \rbrace)\rightarrow 
(\lbrace 1 \rbrace,\lbrace 2 \rbrace)+ (\lbrace 1,\; 2 \rbrace)$. 
For example, this method describes successfully the representation of  
${\cal U}_{\displaystyle{q}}(sl(3))$ of dimension 15 for $m=3$ ($p_{13}=5$,
$p_{23}=2$, $p_{33}=0$).

\vskip 2 truecm

{\bfw Acknowldgements:}  

I thank Daniel Arnaudon and Amitabha Chakrabarti for important and numerous 
discussions and encouragements during this work.

\eject

\centerline{\bfw References.}

\vskip 1 truecm

\item {1.} B. Abdesselam, D. Arnaudon and A. Chakrabarti. 
{\sl Representations of ${\cal U}_{\displaystyle{q}}(sl(N))$ at roots of 
unity,} q-alg/9504006, to be published.

\item {2.} D. Arnaudon and A. Chakrabarti,
{\sl Periodic and partially periodic representations of $SU(N)_{q}$,}
Commun. Math. Phys. {\bf 139} (1991) 461.

\item {3.} T.D. Palev, N.I. Stailova and J. Van der Jeugt, 
{\sl Finite--dimensional representations of the quantum superalgebra 
${\cal U}_{\displaystyle{q}}(gl(n|m))$ and related q-identities},
Commun. Math. Phys. {\bf 166} (1994) 367.

\item {4.} T.D. Palev, Funct. Anal. Appl. vol. 21, No.3 p. 85 (1987) 
(English transl), Journ. Math. Phys. vol. 29, 2589 (1988), Journ. Math. 
Phys. vol. 30, 1433 (1989).

\item {5.} T.D. Palev and V.N. Tolstoy, {\sl Finite dimensional irreducible 
representations of the quantum superalgebra 
${\cal U}_{\displaystyle{q}}(gl(n|1))$},
Commu. Math. Phys. {\bf 141} (1991) 549.

\item {6.} Nguyen Anh Ky and N.I. Stoilova, {\sl Finite dimensional 
representations of the quantum superalgebra 
${\cal U}_{\displaystyle{q}}(gl(2|2))$ II: nontypical representations 
at generic $q$}, hep-th/9411098. 

\item {7.} \v{C}. Burd\'{\i}k, R. C. King and T. A. Welsh,
  {\sl The explicit construction of irreducible representations of the
    quantum algebras ${\cal U}_{\displaystyle{q}}(sl(N))$}
  Proceedings of the 3rd Wigner Symposium, Oxford, September 1993,
  Classical and quantum systems, Eds. L. L. Boyle and A. I. Solomon,
  World Sci., Singapore (1).

\item {8.} V.G. Drinfeld, {\sl Quantum Groups,} Proc.
  Int. Congress of Mathematicians, Berkeley, California, Vol.
  {\bf 1}, Academic Press, New York (1986), 798.

\item {9.} M. Jimbo, {\sl $q$-difference analogue of
    ${\cal U}_{\displaystyle{q}}({\cal G}) $ and the Yang Baxter equation,} Lett. Math. Phys. {\bf 10}
  (1985)  63.

\item {10.} I.M. Gel\'fand and M.L. Testlin. Dokl. Akad. Nauk SSSR, 71, No. 5,
825--828 (1950). 

\item {11.} D. Arnaudon, {\sl Periodic and flat irreducible representations of 
$SU(3)_{q}$}, Comm. Math. Phys. {\bf 134} (1990) 523.

\item {12.} D. Arnaudon and A. Chakrabarti, {\sl Flat periodic 
representations of ${\cal U}_{\displaystyle{q}}({\cal G})$}, 
Commun. Math. Phys. {\bf 139} (1991) 605.

\item {13.} V.K. Dobrev, in Proc. {\sl Int. Group Theory
Conference,} St Andrews, 1989, Vol. 1, Campbell and Robertson (eds.),
London Math. Soc. Lect. Notes Series 159, Cambridge University Press, 1991. 

\item {14.} V.K. Dobrev, {\sl Representations of
    quantum groups,} Proc. ``Symmetries in Science
  V: Algebraic Structures, their Representations, Realizations and
  Physical Applications'', Schloss Hofen, Austria,
  (1990), Eds. B.  Gruber, L.C. Biedenharn
  and H.-D. Doebner (Plenum Press, NY, 1991) pp. 93-135.

\item {15.} G. Lusztig, {\sl Quantum groups at roots of 1}, Geom. Ded. {\bf 35} (1990) 89.


\end